\documentclass[11pt,a4paper]{article}
\usepackage{amsmath,amsfonts,amssymb,amsthm}

\title{Approximation properties of random polytopes\\ associated with Poisson hyperplane processes}
\author{Daniel Hug and Rolf Schneider}
\date{}
\newcommand{\ur}{\mbox{\boldmath$u$}}
\newcommand{\vr}{\mbox{\boldmath$v$}}
\newcommand{\w}{\mbox{\boldmath$w$}}
\newcommand{\x}{\mbox{\boldmath$x$}}
\newcommand{\e}{\mbox{\boldmath$e$}}
\newcommand{\tr}{\mbox{\boldmath$t$}}
\newcommand{\y}{\mbox{\boldmath$y$}}
\newcommand{\z}{\mbox{\boldmath$z$}}
\newcommand{\us}{\mbox{\boldmath$\scriptstyle u$}}

\newcommand{\xs}{\mbox{\boldmath$\scriptstyle x$}}
\newcommand{\ys}{\mbox{\boldmath$\scriptstyle y$}}

\newcommand{\p}{\mbox{\boldmath$p$}}
\newcommand{\q}{\mbox{\boldmath$q$}}
\newcommand{\0}{\mbox{\boldmath$o$}}
\newcommand{\ko}{\mbox{\boldmath$\scriptstyle o$}}
\newcommand{\R}{\mathbb{R}}
\newcommand{\Sd}{\mathbb{S}^{d-1}}
\newcommand{\I}{{\bf 1}}
\newcommand{\fed}{\,\rule{.1mm}{.26cm}\rule{.24cm}{.1mm}\,}

\newcommand{\rP}{{\rm P}}
\newcommand{\rE}{{\rm E} }
\newcommand{\D}{{\rm d}}

\begin{document}
\maketitle

\begin{quote}
{\small
\begin{center}
{\bf Abstract} 
\end{center}

We consider a stationary Poisson hyperplane process with given directional distribution and intensity in $d$-dimensional Euclidean space. Generalizing the zero cell of such a process, we fix a convex body $K$ and consider the intersection of all closed halfspaces bounded by hyperplanes of the process and containing $K$. We study how well these random polytopes approximate $K$ (measured by the Hausdorff distance) if the intensity increases, and how this approximation depends on the directional distribution in relation to properties of $K$.\\[2mm]
}
\end{quote}

\section{Introduction}

Asymptotic properties of the convex hull of $n$ independent, identically distributed random points in $\R^d$, as $n$ tends to infinity, are an actively studied topic of stochastic geometry; see, for example, Subsection 8.2.4 of the book \cite{SW08} and the more recent survey by Reitzner \cite{Rei10}. Very often, one studies uniform random points in a given convex body and measures the rate of approximation by the volume difference, or the difference of other global functionals, or one investigates the asymptotic behaviour of combinatorial quantities such as face numbers. In contrast, approximation by random polytopes, measured in terms of the Hausdorff metric $\delta$, has been investigated less frequently. We recall that the Hausdorff distance of two nonempty compact  sets $K,L\subset {\mathbb R}^d$ is defined by
$$ \delta(K,L)=\max\left\{\max_{\xs\in K}\min_{\ys\in L} \|\x-\y\|,\max_{\xs\in L}\min_{\ys\in K} \|\x-\y\|\right\}.$$
For results on Hausdorff distances of random polytopes we refer to Note 5 for Subsection 8.2.4 in \cite{SW08} and mention here only the following. For a convex body $K$ of class $C^2_+$ (that is, with a twice continuously differentiable boundary with positive Gauss curvature), B\'ar\'any \cite{Bar89} (Theorem 6) showed that the Hausdorff distance from $K$ to the convex hull $K_n$ of $n$ i.i.d. uniform random points in $K$ satisfies
$$ \rE\,\delta(K,K_n)\sim\left(\frac{\log n}{n}\right)^{2/(d+1)}$$
as $n\to\infty$ (here $f(n)\sim g(n)$ means that there are constants $c_1,c_2$ such that $c_1 g(n)<f(n)< c_2g(n)$). A result of D\"{u}mbgen and Walther \cite{DW96} (Corollary 1) says that, for an arbitrary convex body $K$,
$$ \delta(K,K_n)= {\rm O}\left(\left(\frac{\log n}{n}\right)^{1/d}\right)\quad\mbox{almost surely}.$$

The second standard approach to convex polytopes, generating them as intersections of closed halfspaces instead of convex hulls of points, was, for the case of random polygons in the plane, already considered in the third of the seminal papers by R\'{e}nyi and Sulanke \cite{RS63, RS64, RS68}, which initiated this subject. Nevertheless, this approach has later not found equal attention in the study of random polytopes. About the role that duality, either in an exact or a heuristic sense, can play here, we refer to the introduction of \cite{BS10}. This alternative approach has to offer some new aspects, in particular since random hyperplanes naturally come with some directional distribution, which influences the random polytopes that they generate. This aspect is emphasized in the present article, where we consider random polytopes generated by a stationary Poisson hyperplane process, with an arbitrary directional distribution. 

Let $X$ be a stationary nondegenerate (see \cite[p.\,486]{SW08}) Poisson hyperplane process in Euclidean space $\R^d$, $d\ge 2$ (with scalar product $\langle\cdot,\cdot\rangle$ and norm $\|\cdot\|$). The reader is referred to Chapters 3 and 4 of \cite{SW08} for an introduction, and also for some notational conventions used here. In particular, we recall the convention that a simple point process $X$, which is by definition a simple random counting measure, is often identified with its support, which is a locally finite random set.

For a hyperplane $H$ in $\R^d$, not passing through the origin $\0$, we denote by $H^-_{\ko}$ the closed halfspace bounded by $H$ that contains $\0$. The random polytope
$$ Z_0 := \bigcap_{H\in X} H^-_{\ko}$$
is called the {\em zero cell} of $X$ (it is also known as the {\em Crofton polytope} of $X$).

A generalization of this notion is obtained as follows. Let $K\subset \R^d$ be a convex body, by which we understand, in the following, a compact convex subset with interior points. For a hyperplane $H$ not intersecting $K$ we denote by $H^-_K$ the closed halfspace bounded by $H$ that contains $K$. Then we define the $K${\em-cell} of $X$ as the random polytope
$$ Z_K:= \bigcap_{H\in X,\,H\cap K=\emptyset} H^-_K.$$
The almost sure boundedness of $Z_K$ follows as in the proof of \cite[Theorem 10.3.2]{SW08}. 
In the following we are interested in the question how well $K$ is approximated by $Z_K$, if the intensity of the process $X$ tends to infinity. Since the intensity is a constant multiple of the expected number of hyperplanes in the process that hit $K$, the analogy to convex hulls of an increasing number of points is evident.

We consider approximation in sense of the Hausdorff metric $\delta$ on the space ${\mathcal K}^d$ of convex bodies in $\R^d$. Of course, in order that approximation of $K$ by $Z_K$ be possible at all, the convex body $K$ and the directional distribution of the hyperplane process $X$ must somehow be adapted to each other. For example, a ball $K$ cannot be approximated arbitrarily closely by $Z_K$ if the hyperplane process $X$ has only hyperplanes of finitely many directions. To make this more precise, let $N$ be a closed subset of the unit sphere $\Sd$, not contained in a closed halfsphere. For a given convex body $K$, we denote by ${\mathcal P}(K,N)$ the set of all polytopes which are finite intersections of closed halfspaces containing $K$ and with outer unit normal vectors in $N$.

\vspace{3mm}

\noindent{\bf Proposition 1.} {\em The convex body $K$ can be approximated arbitrarily closely, with respect to the Hausdorff metric, by polytopes from ${\mathcal P}(K,N)$ if and only if ${\rm supp}\,S_{d-1}(K,\cdot)\subset N$.}

\vspace{3mm}

Here {\rm supp} denotes the support of a measure, and $S_{d-1}(K,\cdot)$ is the surface area measure of $K$ (see \cite{Sch14}, Section 4.2, for example). We shall give a proof of Proposition 1 in the next section. It serves here only to motivate the assumption (\ref{n1}) made below.

The intensity measure $\Theta= \rE X(\cdot)$ of $X$ is assumed, as usual, to be locally finite. It can then be represented in the form (see \cite{SW08}, (4.33))
\begin{equation}\label{2.0} 
\Theta(A) = 2\gamma\int_{{\mathbb S}^{d-1}}\int_0^\infty \I_A(H(\ur,t))\,\D t\,\varphi(\D \ur)
\end{equation}
for $A\in{\cal B}({\mathcal H}^d)$, where $\gamma>0$ is the intensity and $\varphi$ is the spherical directional distribution of $X$; the latter is an even Borel probability measure on the unit sphere $\Sd$ 
which is not concentrated on a great subsphere. Later, when $\varphi$ is fixed and $\gamma$ varies, we write $\Theta_\gamma$ instead of $\Theta$. By ${\mathcal H}^d$ we denote the space of hyperplanes in $\R^d$, and ${\mathcal B}(T)$ is the $\sigma$-algebra of Borel sets of a topological space $T$. Further,
$$ 
H(\ur,t)=\{\x\in\R^d: \langle\x,\ur\rangle=t\}
$$
for $\ur\in {\mathbb S}^{d-1}$ and $t>0$ is the standard parametrization of a hyperplane not passing through the origin $\0$. For convenience (in view of some later estimations of constants), we also assume that $\gamma\ge 1$.

For $K\in{\mathcal K}^d$, the Hausdorff distance $\delta(K,P)$ of $K$ from a polytope $P$ containing it is the smallest number $\varepsilon\ge 0$ such that $P\subset K(\varepsilon)$, where $K(\varepsilon)=K+\varepsilon B^d$ ($B^d$ is the unit ball) denotes the outer parallel body of $K$ at distance $\varepsilon$. Thus, for given $\varepsilon>0$ the probability $\rP\{\delta(K,Z_K)>\varepsilon \}$, in which we are interested, is equal to $\rP\{Z_K\not\subset K(\varepsilon)\}$. First we give a necessary and sufficient condition that this probability tends to zero if the intensity of the process $X$ tends to infinity; if the condition is satisfied, we obtain that the decay is exponential. Under a slightly stronger assumption, this can then be used to derive our main results, concerning the rate of convergence.

We assume in the following that the surface area measure of the given convex body $K$ satisfies
\begin{equation}\label{n1}
{\rm supp}\,S_{d-1}(K,\cdot) \subset {\rm supp}\,\varphi.
\end{equation}
By Proposition 1, this assumption is necessary for arbitrarily good approximation of $K$ by $Z_K$. Theorem 1 shows, in a stronger form, that it is also sufficient.

For $\y\in \R^d\setminus K$, let $K^{\ys}:= {\rm conv}(K\cup\{\y\})$. For $\varepsilon>0$ we define
\begin{equation}\label{n4} 
\mu(K,\varphi,\varepsilon):= \min_{\ys\in {\rm bd}\,K(\varepsilon)} \int_{\Sd}[h(K^{\ys},\ur)-h(K,\ur)]\,\varphi(\D \ur),
\end{equation}
where $h$ denotes the support function. Lemma 1, to be proved in the next section, shows that condition (\ref{n1}) implies $\mu(K,\varphi,\varepsilon)>0$.

\vspace{3mm}

\noindent{\bf Theorem 1.} {\em Let $K\in{\mathcal K}^d$ be a convex body. Let $X$ be a stationary Poisson hyperplane process in ${\mathbb R}^d$ with intensity $\gamma$ and with a directional distribution $\varphi$ satisfying $(\ref{n1})$. There are positive constants $C_1(\varepsilon), C_2$ (both depending  on $K$, $\varphi$, $d$) such that the following holds. If $0<\varepsilon\le 1$, then}
\begin{equation}\label{5.0} 
\rP\left\{\delta(K,Z_K)>\varepsilon \right\} \le C_1(\varepsilon)\exp\left[-C_2 \mu(K,\varphi,\varepsilon)\gamma\right],
\end{equation}
where $\mu(K,\varphi,\varepsilon)>0$.

\vspace{3mm}

In order to be able to deal with convergence for increasing intensities, we consider an embedding of the stationary Poisson hyperplane processes $X_\gamma$ with intensity $\gamma>0$, directional distribution $\varphi$ and intensity measure 
$$ \rE X_\gamma(\cdot)  = 2\gamma \int_{{\mathbb S}^{d-1}}\int_0^\infty \I\{H(\ur,t)\in\cdot\}\,\D t\,\varphi(\D \ur)=: \Theta_\gamma$$ 
into a Poisson process $\xi$ in $[0,\infty)\times\mathcal{H}^d$ (on a suitable probability space) with intensity measure $\lambda\otimes \Theta_1$, where $\lambda$ denotes Lebesgue measure on $[0,\infty)$. Then $\xi([0,\gamma]\times\cdot)$ is a Poisson hyperplane process in ${\mathbb R}^d$ with intensity measure $\Theta_\gamma$, thus $X_\gamma$ is stochastically equivalent to $\xi([0,\gamma]\times\cdot)$ (e.g., by \cite{SW08}, Theorem 3.2.1). In the following, we can identify $X_\gamma$ with $\xi([0,\gamma]\times\cdot)$. Let $Z_K^{(\gamma)}$ denote the $K$-cell associated with $\xi([0,\gamma]\times\cdot)$. Then we have $K\subset Z_K^{(\tau)}\subset Z_K^{(\gamma)}$ for $\tau\ge \gamma>0$, and therefore $\delta(K,Z_K^{(\tau)})\le \delta(K,Z_K^{(\gamma)})$. This shows that 
$$ \rP\left\{\sup_{\tau\ge \gamma}\delta(K,Z_K^{(\tau)})\ge \varepsilon\right\}=\rP\left\{\delta(K,Z_K^{(\gamma)})\ge \varepsilon\right\} \le C_1(\varepsilon)\exp\left[-C_2 \mu(K,\varphi,\varepsilon)\gamma\right] $$
for all $\varepsilon>0$, and thus 
$$ \lim_{\gamma\to \infty} \delta(K,Z_K^{(\gamma)}) = 0 $$
holds almost surely. We state this as a corollary.

\vspace{3mm}

\noindent{\bf Corollary.} {\em If the Poisson hyperplane processes $X_\gamma$, $\gamma\ge 1$, are defined as above on a common probability space and if $Z_K^{(\gamma)}$ denotes the $K$-cell of $X_\gamma$ for a convex body $K\in {\mathcal K}^d$, then condition $\rm (\ref{n1})$ is necessary and sufficient in order that}
\begin{equation}\label{n6} 
\lim_{\gamma\to \infty} \delta(K,Z_K^{(\gamma)}) = 0 \quad\mbox{almost surely}.
\end{equation}

\vspace{3mm}

In the following, we will be interested in rates of convergence. For this, we consider the sequence $X_1,X_2,\dots$ of Poisson hyperplane processes defined as above, with spherical directional distribution $\varphi$, where $X_n$ has intensity $n$. 

Under the sole assumption (\ref{n1}), no statement stronger than (\ref{n6}), involving also a rate of convergence, is possible. In fact, if any decreasing sequence $(\varepsilon_n)_{n\in{\mathbb N}}$ with $\varepsilon_n\to 0$ for $n\to\infty$ is given and if $K$ is a convex body, then the directional distribution $\varphi$ of the hyperplane processes $X_n$ can be chosen in such a way that (\ref{n1}) is satisfied but
\begin{equation}\label{C1}
\rP\{\delta(K,Z_K^{(n)})\ge\varepsilon_n \mbox{ for almost all } n\}=1.
\end{equation}
We prove this at the end of the paper. 

Therefore, no assumption on the convex body $K$ alone allows us to estimate the rate of convergence of $\delta(K,Z_K^{(n)})$ for arbitrary directional distributions $\varphi$. On the other hand, suitable assumptions on the directional distribution, for example
\begin{equation}\label{n7}
\varphi \ge b\sigma
\end{equation}
with a constant $b>0$, where $\sigma$ denotes spherical Lebesgue measure, permit to estimate the rate of convergence for arbitrary convex bodies. This is shown by the first assertion of Theorem 2. 

If the directional distribution does not satisfy such a strong assumption, then rates of convergence can only be estimated if this distribution is suitably adapted to the given convex body. In this sense, we assume that
\begin{equation}\label{n3}
\varphi \ge bS_{d-1}(K,\cdot)
\end{equation}
with some constant $b$. 

If $(Y_n)_{n\in\mathbb{N}}$ is a sequence of real random variables and $f(n)_{n\in\mathbb{N}}$ is a sequence of nonnegative real numbers, we write $Y_n={\rm O}(f(n))$ {\em almost surely} if there is a constant $C<\infty$ such that with probability one we have $Y_n\le Cf(n)$ for sufficiently large $n$. Moreover, we write $Y_n\sim f(n)$ {\em almost surely} if there are constants $0<c\le C<\infty$ such that  with probability one we have $cf(n)\le Y_n\le Cf(n)$ for all sufficiently large $n$. A `ball' in the following is a Euclidean ball of positive radius. One says that a convex body $M$ {\em slides freely} inside a convex body $K$ if $K$ is the union of all translates of $M$ that are contained in $K$.

\vspace{3mm}

\noindent{\bf Theorem 2.} {\em Let $K\in{\mathcal K}^d$ be a convex body. Let $X$ be a stationary Poisson hyperplane process in ${\mathbb R}^d$ with intensity $\gamma$ and with a directional distribution $\varphi$ satisfying $\rm (\ref{n7})$ or $\rm (\ref{n3})$. Then
\begin{equation}\label{5.0z}
\delta(K,Z_K^{(n)})={\rm O}\left(\left(\frac{\log n}{n}\right)^{1/d}\right)\quad\mbox{almost surely,}
\end{equation}
as $n\to\infty$. 

Suppose that $\rm (\ref{n3})$ holds. If a ball slides freely inside $K$, then the exponent $1/d$ in $(\ref{5.0z})$ can be replaced by $2/(d+1)$, and if $K$ is a polytope, then it can be replaced by $1$.}

\vspace{3mm}

Under stronger assumptions on $K$ and $\varphi$, we can determine the exact asymptotic order of approximation. 

\vspace{3mm}

\noindent{\bf Theorem 3.} {\em Let the convex body $K\in{\mathcal K}^d$ be such that a ball slides freely inside $K$ and that $K$ slides freely inside a ball. Suppose that the directional distribution $\varphi$ of the stationary Poisson hyperplane processes $X_n$ satisfies 
\begin{equation}\label{n9}
a\sigma \ge \varphi\ge b\sigma
\end{equation} 
with some positive constants $a,b$. Then
\begin{equation}\label{n8}
\delta(K,Z_K^{(n)})\sim \left(\frac{\log n}{n}\right)^{2/(d+1)}\quad\mbox{almost surely,}
\end{equation}
as $n\to\infty$. }

\vspace{3mm}

Note that Theorem 3 covers, in particular, the case where $K$ is of class $C^2_+$ and the hyperplane processes $X_n$ are isotropic, that is, their directional distribution $\varphi$ is invariant under rotations and thus is equal to the normalized spherical Lebesgue measure. If $K$ is of class $C^2_+$, then the assumptions on $K$ are satisfied by Blaschke's rolling theorem (Corollary 3.2.13 in \cite{Sch14}).

In the next section, we prove some auxiliary results. Theorem 1 is proved in Section 3, and the proofs of Theorems 2 and 3 follow in Section 4.

\section{Auxiliary results}

\noindent{\em Proof of Proposition $1$.} By \cite{Sch14}, Theorem 4.5.3, the support of the area measure $S_{d-1}(K,\cdot)$ is equal to ${\rm cl\;extn}\, K$, the closure of the set of extreme (unit) normal vectors of $K$. 

Suppose now that $K$ can be approximated arbitrarily closely by polytopes from ${\mathcal P}(K,N)$. Let $\x$ be a regular boundary point of $K$, and let $(\x_i)_{i\in{\mathbb N}}$ be a sequence of points in $\R^d\setminus K$ converging to $\x$. To each $i$, there exists a polytope $P_i\in {\mathcal P}(K,N)$ not containing $\x_i$, hence there is a closed halfspace $H^-_i$ with outer normal vector $\ur_i\in N$ containing $K$ but not $\x_i$. For $i\to\infty$, the sequence of hyperplanes $H_i$ bounding $H^-_i$ has a convergent subsequence; its limit is the unique supporting hyperplane of $K$ at $\x$. It follows that the outer unit normal vector of $K$ at $\x$ belongs to the closed set $N$. A normal vector at a regular boundary point of $K$ is a $0$-exposed normal vector. Since $\x$ was an arbitrary regular boundary point of $K$, the set $N$ contains the set of $0$-exposed normal vectors of $K$. The closure of the $0$-exposed normal vectors is equal to the closure of the extreme normal vectors (see Theorem 2.2.9 of \cite{Sch14}, also for the terminology used here). Hence, ${\rm cl\;extn}\, K\subset N$.

Conversely, suppose that ${\rm cl\;extn}\, K\subset N$. The body $K$ is the intersection of its supporting halfspaces with a regular point of $K$ in the boundary (see \cite{Sch14}, Theorem 2.2.5). The outer unit normal vector of such a halfspace is extreme and hence belongs to $N$. Thus, denoting by $H^-(K,\ur)$ the supporting halfspace of $K$ with outer unit normal vector $\ur$, we have $K=\bigcap_{\us\in N}H^-(K,\ur)$. Therefore, if $\varepsilon>0$, then
$$ \bigcap_{\us\in N} {\rm bd}(K+\varepsilon B^d)\cap H^-(K,\ur)=\emptyset.$$
By compactness, there is a finite subset $F\subset N$ such that the corresponding intersection is empty, which implies that
$$ P:= \bigcap_{\us\in F}H^-(K,\ur)\subset{\rm int}(K+\varepsilon B^d ).$$
Thus, $P$ is a polytope in ${\mathcal P}(K,N)$ with $\delta(K,P)<\varepsilon$. Since $\varepsilon>0$ was arbitrary, this shows that $K$ can be approximated arbitrarily closely by polytopes from ${\mathcal P}(K,N)$.
\qed

\vspace{3mm}

In the rest of this paper, $c_1,c_2,\dots$ denote positive constants that depend only on $K$, $\varphi$ and the dimension $d$.

\vspace{3mm}

\noindent{\bf Lemma 1.} {\em Let $K\in{\mathcal K}^d$ and let $\varphi$ be a probability measure on $\Sd$. Let $0<\varepsilon\le 1$. 

\noindent $\rm (a)$ If $(\ref{n1})$ holds, then $ \mu(K,\varphi,\varepsilon)>0$.

\noindent $\rm (b)$ If $\rm (\ref{n7})$ holds, then there exists a constant $c_1$ such that
\begin{equation}\label{5.0a} 
\mu(K,\varphi,\varepsilon) \ge c_1\varepsilon^d.
\end{equation}

In $\rm (c),\,(d),\,(e)$ it is assumed that $(\ref{n3})$ is satisfied. 

\noindent $\rm (c)$ For $\varepsilon\le D(K)$, where $D(K)$ denotes the diameter of $K$, there exists a constant $c_2$ such that
\begin{equation}\label{5.1} 
\mu(K,\varphi,\varepsilon) \ge c_2\varepsilon^d.
\end{equation}

\noindent $\rm (d)$ If a ball slides freely inside $K$, then there exists a constant $c_3$ such that
\begin{equation}\label{5.2}
\mu(K,\varphi,\varepsilon)\ge c_3\varepsilon^{(d+1)/2}.
\end{equation}

\noindent $\rm (e)$ If $K$ is a polytope, then there exists a constant $c_4$ such that
\begin{equation}\label{5.3}
\mu(K,\varphi,\varepsilon) \ge c_4\varepsilon.
\end{equation} 
}

\noindent{\em Proof.} 
(a) Let $(\ref{n1})$ be satisfied. Let $\y\in\R^d\setminus K$. Let $V_d$ denote the volume and $V$ the mixed volume in $\R^d$. Using a formula for mixed volumes (\cite{Sch14}, (5.19)) and Minkowski's inequality (e.g., \cite{Sch14}, (7.18)), we get
\begin{eqnarray*}
& & \frac{1}{d} \int_{\Sd} [h(K^{\ys},\ur)-h(K,\ur)]\,S_{d-1}(K,\D\ur)\\
&& =V(K^{\ys},K,\dots,K)-V_d(K)\\
&&\ge V_d(K^{\ys})^{\frac{1}{d}}V_d(K)^{\frac{d-1}{d}}-V_d(K)\\
&&=V_d(K)^{\frac{d-1}{d}}\left[V_d(K^{\ys})^{\frac{1}{d}}-V_d(K)^{\frac{1}{d}}\right]\\
&&>0.
\end{eqnarray*}
The integrand is nonnegative and continuous as a function of $\ur$. Since the integral is positive, there exists a neighbourhood (in $\Sd$) of some point $\ur_0\in{\rm supp}\,S_{d-1}(K,\cdot)$ on which the integrand is positive. By (\ref{n1}), $\ur_0\in {\rm supp}\,\varphi$, and hence
$$ g(\y):=\int_{\Sd} [h(K^{\ys},\ur)-h(K,\ur)]\,\varphi(\D\ur)>0.$$
The function $g$ is continuous, hence on each compact subset of $\R^d\setminus K$ it attains a minimum. This proves that $\mu(K,\varphi,\varepsilon)>0$.

(b) Suppose that (\ref{n7}) holds. For the proof of (\ref{5.0a}), let $K\in{\mathcal K}^d$ be given. Let $\y\in {\rm bd}\, K(\varepsilon)$ and let $\x$ be the point in $K$ nearest to $\y$. Then $N(\x):=(\y-\x)/\varepsilon$ 
is an outer unit normal vector of $K$ at $\x$. We denote by $H^-$ the closed halfspace bounded by the hyperplane through $\x$ and orthogonal to $N(\x)$ and containing $K$. If $D(K)$ denotes the diameter of $K$, then $K\subset H^-\cap (\x+D(K)B^d)$. Define $\beta=\beta(\varepsilon)\in [0,\pi/2)$ by $\cos\beta=D(K)/\sqrt{D(K)^2+\varepsilon^2}$ and let $S(\y,\varepsilon)$ be the set of all ${\ur}\in\Sd$ such that $\angle ({\ur},N(\x) )\le\beta/2$. Then 
\begin{equation}\label{eqle1}
\sigma(S(\y,\varepsilon))\ge c_5\sin^{d-1}(\beta/2)\ge c_6\varepsilon^{d-1}.
\end{equation}
For $\ur\in S(\y,\varepsilon)\setminus \{N(\x)\}$ there is a unique unit vector $\e$ orthogonal to $N(\x)$ such that $\ur=\tau N(\x)+\sqrt{1-\tau^2}\,\e$ with $0<\tau<1$. 
With $\z:=D(K)\e$ we then obtain 
\begin{align}
h(K^{\ys},\ur)-h(K,\ur)&\ge\langle \y,\ur\rangle-\langle \z,\ur\rangle =\langle \y-\z,\ur\rangle\nonumber\\
&\ge D(K)\left\langle \frac{\y-\z}{\|\y-\z\|},\ur\right\rangle\ge D(K)\sin(\beta/2)\nonumber\\
&\ge c_7\varepsilon,\label{le1eq2}
\end{align}
for all $\ur\in S(\y,\varepsilon)$. 
Combining (\ref{n7}), \eqref{eqle1} and \eqref{le1eq2}, we obtain 
\begin{eqnarray*}
& &  \frac{1}{b}\int_{\Sd}[h(K^{\ys},\ur)-h(K,\ur)]\, \varphi(\D\ur)\\
& & \ge\int_{\Sd}[h(K^{\ys},\ur)-h(K,\ur)]\, \sigma(\D\ur)\ge \sigma(S(\y,\varepsilon)) c_7\varepsilon\ge c_8\varepsilon^d,
\end{eqnarray*}
which completes the proof of (b).

Now suppose that $(\ref{n3})$ holds. From the estimate in the proof of (a) we get
\begin{eqnarray*}
\frac{1}{bd} \int_{\Sd} [h(K^{\ys},\ur)-h(K,\ur)]\,\varphi(\D\ur)
&\ge& \frac{1}{d} \int_{\Sd} [h(K^{\ys},\ur)-h(K,\ur)]\,S_{d-1}(K,\D\ur)\\
& \ge&V_d(K)^{\frac{d-1}{d}}\left[V_d(K^{\ys})^{\frac{1}{d}}-V_d(K)^{\frac{1}{d}}\right]\\
&\ge& c_9\left[V_d(K^{\ys})-V_d(K)\right].
\end{eqnarray*}

(c) For the proof of (\ref{5.1}), let $\y\in{\rm bd}\,K(\varepsilon)$ and let $C$ be the cone with apex $\y$ spanned by $K$. Let $\y'$ be the point in $K$ nearest to $\y$. The vector $\y-\y'$ has length $\varepsilon$, and the hyperplane $H'$ orthogonal to it and passing through $\y'$ supports $K$. Let $H$ be the other supporting hyperplane of $K$ parallel to $H'$. Let $\Delta$ be the convex hull of $\y$ and $H\cap C$ and $\Delta'$ the convex hull of $\y$ and $H'\cap C$. Denoting by $D(K)$ the diameter of $K$ and assuming that $\varepsilon\le D(K)$, we have
$$ V_d(K^{\ys})-V_d(K) \ge V_d(\Delta') \ge \left(\frac{\varepsilon}{D(K)+\varepsilon}\right)^d V_d(\Delta) \ge  \left(\frac{\varepsilon}{2D(K)}\right)^dV_d(K).$$
This gives (\ref{5.1}).

(d) Suppose that a ball of radius $r>0$ slides freely inside $K$. Since $\mu(\cdot,\varphi,\varepsilon)$ is translation invariant, we can assume that $K$ contains the ball $B(\0,r)$ of radius $r$ centred at $\0$. Let 
$R>0$ be such that $K\subset B(\0,R)$. For $s>0$, the convex body
$$ K^s:=\{\x\in\R^d:V_d(K^{\xs})-V_d(K)\le s\} $$
is known as an illumination body of $K$ (cf.~\cite[p.~258]{Werner1994}; the convexity follows from Satz 4 in F\'ary and R\'edei \cite{FR50}). Now let $\y\in \text{bd}\, K(\varepsilon)$ and put $\nu:=V_d(K^{\ys})-V_d(K)$, then $\y\in\text{bd}\, K^\nu$. Let $\x\in\text{bd}\, K$ be determined by $\{\x\}=[\0,\y]\cap\text{bd}\,K$, and denote by $N(\x)$ the unique exterior unit normal vector of $K$ at $\x$ (the normal vector is unique since by assumption there is a ball $B'$ of radius $r>0$ with $\x\in B'\subset K$). Since $B(\0,r)\subset K$, we have
$$ \langle \x,N(\x)\rangle\ge r,\qquad \langle \x/\|\x\|,N(\x)\rangle\ge r/R. $$
From $\|\y\|-\|\x\| \ge\varepsilon$ we get $\|\y\|^d-\|\x\|^d \ge dr^{d-1}\varepsilon$. Therefore, Lemma 2 in \cite{Werner1994} yields
$$ \nu^{2/(d+1)} \ge c_{10} rr^{(d-1)/(d+1)}\left(\left(\frac{\|\y\|}{\|\x\|}\right)^d-1\right)
\ge c_{11}R^{-d}\left(\|\y\|^d-\|\x\|^d\right) \ge c_{12}\varepsilon,$$
hence
$$
V_d(K^{\ys})-V_d(K)\ge c_{13}\,\varepsilon^{(d+1)/2},
$$
which gives (\ref{5.2}).

(e) Now suppose that $K$ is a polytope. Let $\y\in{\rm bd}\,K(\varepsilon)$ and let $\y'$ be the point in $K$ nearest to $\y$. Put $\vr:=(\y-\y')/\|\y-\y'\|$, and let $F$ denote the unique (proper) face of $K$ which contains $\y'$ in its relative interior. Let $F_1,\dots,F_m$ be the facets of $K$ that contain $F$, and let $\ur_1,\dots,\ur_m$ be their outer unit normal vectors. By \cite[p. 85 and Theorem 2.4.9]{Sch14}, we have
$$
\vr\in N(K,F)=N(K,\y')=\text{pos}\{\ur_i:i=1,\ldots,m\},
$$
where $N(K,F)$ and $N(K,\y')$ are the normal cones of $K$ at $F$ and $\y'$, respectively, and $\text{pos}$ denotes the positive hull. For any unit vector $\w\in N(K,F)$ there is some $i\in\{1,\ldots,m\}$ such that $\langle \w, \ur_i\rangle>0$; in particular, 
$$
a(F,\w):=\max\{\langle \w,\ur_i\rangle:i=1,\dots,m\}>0
$$ 
and $a(F,\vr)=\langle \vr,\ur_{i_0}\rangle>0$ for some $i_0\in\{1,\dots,m\}$. Since $N(K,F)\cap \Sd$ is compact, we have
$$
a(F):=\min\{a(F,\w):\w\in N(K,F)\cap \Sd\}>0
$$
and thus
$$
c_{14}:=\min\{a(F):F\text{ is a proper face of }K\}>0.
$$
Therefore, with $c_{15}:=\min\{V_{d-1}(F):F\text{ is a facet of }K\}>0$, where $V_{d-1}$ denotes the $(d-1)$-dimensional volume, we get
\begin{align*}
\int_{\Sd} [h(K^{\ys},\ur)-h(K,\ur)]\,S_{d-1}(K,\D\ur)\ge&  \langle \y-\y', \ur_{i_0}\rangle V_{d-1}(F_{i_0})\\
\ge& \|\y-\y'\|\cdot c_{14} c_{15}=c_{16}\varepsilon .
\end{align*}
This yields (\ref{5.3}).\qed

\vspace{3mm}

\noindent{\bf Remark.} Although in the case of a general convex body $K$, the derivation of the estimate (\ref{5.1}) may seem rather crude, the order of $\varepsilon^d$ cannot be improved. In fact, if (\ref{5.1}) would be replaced by $\mu(K,\varphi,\varepsilon) \ge c_2\varepsilon^\alpha$ with $1<\alpha< d$, then a counterexample would be provided by a body $K$ which in a neighbourhood of some boundary point is congruent to a suitable part of a body of revolution with meridian curve given by $\mu(t)= |t|^r$ with $1< r< \frac{d-1}{\alpha-1}$.

\vspace{3mm}
\noindent{\bf Lemma 2.} {\em Let the convex body $K\in{\mathcal K}^d$ be such that a ball slides freely inside $K$. Assume further that
\begin{equation}\label{n9a}
a\sigma\ge\varphi
\end{equation} 
with some positive constant $a$. Then
$$  \int_{\Sd} [h(K^{\ys},\ur)-h(K,\ur)]\,\varphi(\D\ur)\le c_{17} \varepsilon^{(d+1)/2} $$
for $\varepsilon>0$ and $\y\in{\rm bd}\,K(\varepsilon)$.}

\vspace{3mm} 

\noindent{\em Proof.} Let $\y\in {\rm bd}\,K(\varepsilon)$. From (\ref{n9a}) we get
$$ \int_{\Sd}[h(K^{\ys},\ur)-h(K,\ur)]\,\varphi(\D \ur) \le c_{18}\int_{\Sd}[h(K^{\ys},\ur)-h(K,\ur)]\,\sigma(\D \ur). $$

Let $\x$ be the point in $K$ nearest to $\y$; then $\y=\x+\varepsilon N(\x)$, where $N(\x)$ is the outer unit normal vector of $K$ at $\x$. By assumption, a ball, say of radius $r>0$, slides freely inside $K$. In particular, some ball $B$ of radius $r$ satisfies $\x\in B\subset K$. Let
$$ {\rm Cap}\,(\y,\varepsilon):=\left\{\ur\in\Sd:\langle \ur,N(\x) \rangle\ge \frac{r}{r+\varepsilon}\right\}.$$
For $\ur\in {\mathbb S}^{d-1}\setminus {\rm Cap}\,(\y,\varepsilon)$ we have $h(K^{\ys},\ur) -h(K,\ur)=0$. If  $h(K^{\ys},\ur) -h(K,\ur)\not=0$, then
$$ h(K^{\ys},\ur) -h(K,\ur) \le \langle \y-\x,\ur\rangle \le\varepsilon. $$
With $\alpha(\varepsilon):=\arccos r/(r+\varepsilon)$ this gives
\begin{eqnarray*}
\int_{\Sd}[h(K^{\ys},\ur)-h(K,\ur)]\,\sigma(\D \ur) &\le&
\int_{{\rm Cap}\,(\ys,\varepsilon)}\varepsilon\,\sigma(\D\ur)\\
&\le& c_{18} \varepsilon \sin^{d-1}\alpha(\varepsilon)= c_{18}\varepsilon\sqrt{1-(r/(r+\varepsilon))^2}^{\,d-1}\\
&\le&  c_{19}\varepsilon^{(d+1)/2}.
\end{eqnarray*}
This yields the assertion.\qed

The following lemma is sufficient for our purpose; it does not aim at an optimal order.

\vspace{3mm}

\noindent{\bf Lemma 3.} {\em Let $K\in{\mathcal K}^d$ be a convex body which slides freely in some ball. There are constants $c_{20}, c_{21} >0$ such that the following holds. For $0< \varepsilon< c_{20}$, let $m(\varepsilon)$ be the largest number $m$ such that there are $m$ points in ${\rm bd}\,K(\varepsilon)$ with the property that each segment connecting any two of them intersects the interior of $K$. Then }
$$ m(\varepsilon)\ge c_{21}\varepsilon^{-1/2}.$$

\noindent{\em Proof.} The convex body $K$ (which has interior points, by our general assumption) contains some ball, without loss of generality the ball $rB^d$. Let $R$ be such that $K$ slides freely in a ball of radius $R$. We put $c_{20}:= \min\{2R,(\pi r)^2/64R\}$ and assume that $0< \varepsilon< c_{20}$.

For points $\x,\y\in{\rm bd}\,K(\varepsilon)$, we assert that 
\begin{equation}\label{2.10}
\|\x-\y\|\ge 4\sqrt{R\varepsilon}\enspace \Rightarrow\enspace [\x,\y]\cap{\rm int}\,K\not=\emptyset.
\end{equation}
For the proof, let $\x,\y\in{\rm bd}\,K(\varepsilon)$ and suppose that $ [\x,\y]\cap{\rm int}\,K=\emptyset$. Let $\p\in K$ and $\q\in{\rm aff}\,\{\x,\y\}$ be points of smallest distance. If $\p\not=\q$, then the hyperplane $H$ through $\p$ orthogonal to $\q-\p$ supports $K$. If $\p=\q$, then the line through $\x$ and $\y$ touches $K$, and we choose $H$ as a supporting hyperplane of $K$ containing that line. The body $K$ slides freely in a ball, say $B$, of radius $R$, hence $K$ is a summand of $B$ (\cite{Sch14}, Theorem 3.2.2). This means that there exists a compact convex set $M\subset{\mathbb R}^d$ such that $K+M=B$. 

Let $\ur$ denote the outer unit normal vector of the supporting hyperplane $H$ of $K$ at $\p$, so that $h(K,\ur)=\langle \p,\ur\rangle$. There is a point $\tr\in M$ with $h(M,\ur)=\langle \tr,\ur\rangle$, and the point $\z:= \p+\tr$ satisfies $\z\in B$ and $h(B,\ur)=\langle\z,\ur\rangle$. It follows that $K\subset B-\tr$ and that $H$ is a supporting hyperplane of $B-\tr$ at $\p$. 

The ball $(B-\tr)+\varepsilon B^d$ contains $K(\varepsilon)$ and hence the segment $[\x,\y]$. The line parallel to $[\x,\y]$ through $\p$ lies in $H$ and intersects the ball $(B-\tr)+\varepsilon B^d$ in a segment $S$, which is not shorter than $[\x,\y]$. Thus, $\|\x-\y\|\le{\rm length}(S)= 2\sqrt{2R\varepsilon+\varepsilon^2} <4\sqrt{R\varepsilon}$, since $\varepsilon< 2R$. This proves (\ref{2.10}).
 
Let $m$ be the largest integer with
$$ m\le \frac{\pi r}{4\sqrt{R}}\,\varepsilon^{-1/2}.$$
Then $m\ge 2$ (by the choice of $c_{20}$), and there is a constant $c_{21}$ with $m\ge c_{21}/\sqrt{\varepsilon}$. Let $C$ be an arbitrary great circle of the ball $rB^d$. On $C$, we choose $m$ equidistant points $\y_1,\dots,\y_m$. For $i\not=j$ we have $\|\y_i-\y_j\|\ge 2r\sin(\pi/m)>r\pi/m$. Let $\x_i= \lambda_i\y_i\in{\rm bd}\,K(\varepsilon)$ with $\lambda_i>0$, then $\lambda_i>1$ for $i=1,\dots,m$ and hence $\|\x_i-\x_j\|>r\pi/m\ge 4\sqrt{R\varepsilon}$ for $i\not=j$. By (\ref{2.10}), this completes the proof. \qed

\section{Proof of Theorem 1}

We assume that $X$ and $K$ are as in Theorem 1 and satisfy the assumptions mentioned above, that is, $\varphi$ is not concentrated on a great subsphere, $\gamma\ge 1$, and the inclusion (\ref{n1}) holds. Without loss of generality, we may assume that $\0\in{\rm int}\,K$. Recalling that $Z_0$ denotes the zero cell of $X$, we note that by the independence properties of the Poisson process we have
$$ \rP\{Z_K\not\subset K(\varepsilon)\}= \rP\left\{Z_0\not\subset K(\varepsilon)\mid K\subset Z_0\right\}.$$
The conditional probability involving the zero cell is slightly more convenient to handle. 

For a compact convex set $L\subset\R^d$ we define
$$ {\cal H}_L:=\{H\in {\mathcal H}^d:H\cap L\not=\emptyset\}$$
and 
$$ \Phi(L):= \Theta({\cal H}_L). $$
By (\ref{2.0}) we have 
\begin{equation}\label{n2}
\Phi(L) = 2\gamma\int_{\Sd} h(L,\ur)\,\varphi(\D\ur). 
\end{equation}

The following two lemmas use ideas from the proofs of Lemmas 3 and 5 in \cite{HS07}, but the present situation is simpler. As there, we use the abbreviation
$$ H_1^-\cap\dots\cap H_n^-=:P(H_{(n)}),$$
where $H_1,\dots,H_n$ are hyperplanes not passing through $\0$ and $H_i^-$ is the closed halfspace bounded by $H_i$ that contains $\0$.

Let $\|\x\|_K= \min\{\lambda\ge 0: \x\in\lambda K\}$ for $\x\in {\mathbb R}^d$. For a nonempty compact convex set $L$, we define $\|L\|_K:=\max\{\|\x\|_K:\x\in L\}$. For $\varepsilon\ge 0$ and $m\in {\mathbb N}$, let
$$ {\cal K}^d_\varepsilon(m):=\{L\in{\cal K}^d: K\subset L \not\subset K(\varepsilon),\,\|L\|_K\in (m,m+1]\}$$
and
$$ q_\varepsilon(m):= \rP\{Z_0\in{\cal K}^d_\varepsilon(m)\}.$$
We abbreviate
$$(m+1)K=:K_m.$$
We have
\begin{equation}\label{1} 
q_\varepsilon(m) = \sum_{N=d+1}^\infty \rP\{X({\cal H}_{K_m})=N\}p(N,m,\varepsilon)
\end{equation}
with 
\begin{eqnarray*}
p(N,m,\varepsilon) &:=& \rP\{Z_0\in{\cal K}^d_\varepsilon(m)\mid X({\cal H}_{K_m})=N\}\\
&=& \Phi(K_m)^{-N} \int_{{\cal H}_{K_m}^N}{\bf 1}\{P(H_{(N)})\in{\cal K}^d_\varepsilon(m)\}\, \Theta^N(\D(H_1,\dots,H_N)),
\end{eqnarray*}
the latter by a well-known property of Poisson processes (e.g., \cite[Th. 3.2.2(b)]{SW08}),
and
\begin{equation}\label{2} 
\rP\{X({\cal H}_{K_m})=N\} =  \frac{\Phi(K_m)^{N}}{N!} \exp\left[-\Phi(K_m)\right].
\end{equation}

\vspace{3mm}

\noindent{\bf Lemma 4.} {\em There exists a number $m_0$, depending only on $K$, $\varphi$ and $d$, such that
$$ q_0(m) \le c_{22}\exp[-\Phi(K)-c_{23}\gamma  m] $$
for $m\ge m_0$. }

\vspace{3mm}

\noindent{\em Proof.} We modify and adapt the proof of Lemma 3 in \cite{HS07}. If $H_1,\dots,H_N\in {\cal H}_{K_m}$ and if $P:= P(H_{(N)})\in {\cal K}^d_0(m)$, then $P$ 
has a vertex $\vr$ with $m<\|\vr\|_K\le m+1$. Since $\vr$ is the intersection of some $d$ facets of $P$, there exists a $d$-element set $J\subset \{1,\dots,N\}$ with
$$ \{\vr\}= \bigcap_{j\in J}H_j.$$
We denote the segment $[\0,\vr]$ by $S=S(H_i,\, i\in J)$ (where it is assumed that the hyperplanes $H_i$, $i\in J$, have linearly independent normal vectors) and note that
$$ H_{i} \cap {\rm relint}\, S = \emptyset  \qquad\mbox{for}\qquad i=1,\dots,N. $$
For any segment $S=[\0,\vr]$ with $\|\vr\|_K\ge  m$ we have (writing $a^+:=\max\{a,0\}$)
$$ \Phi(S) = 2\gamma \int_{{\mathbb S}^{d-1}} \langle \vr,\ur\rangle^+ \varphi(\D \ur) \ge 2c_{24}\gamma  m$$
with a positive constant $c_{24}$. This follows from the fact that the function
$$ \vr_1\mapsto \int_{{\mathbb S}^{d-1}} \langle \vr_1, \ur\rangle^+ \varphi(\D \ur),\qquad \vr_1\in{\mathbb S}^{d-1},$$
is positive (since $\varphi$ is not concentrated on a great subsphere) and continuous.  Let $m_0$ be the smallest integer $\ge (2/c_{24})\int_{{\mathbb S}^{d-1}}h(K,\ur)\,\varphi(\D \ur)$. 
For  $m\ge m_0$ we then have 
$$ \Phi(S) \ge \Phi(K) +c_{24}\gamma m,$$
and hence
$$ \int_{{\cal H}_{K_m}} {\bf 1} \{H \cap S = \emptyset\} \,\Theta(\D H) = \Phi(K_m) - \Phi(S)\le \Phi(K_m) - \Phi(K) -c_{24}\gamma m,$$
where we used that $S\subset K_m$, since $\|\vr\|_K\le m+1$. 
Now we obtain
\begin{eqnarray*}
p(N,m,\varepsilon) & \le & \binom{N}{d}\Phi(K_m)^{-N} \int_{{\cal H}_{K_m}^{d}}  {\bf 1} \left\{ \|S(H_j,\,j\in\{1,\dots,d\})\|_K \ge  m\right\} \\
&  & \int_{{\cal H}_{K_m}^{N-d}} {\bf 1} \left\{H_i \cap  S(H_j,\,j\in\{1,\dots,d\})= \emptyset \mbox{ for } i=d+1,\dots,N\right\} \\
&  & \times \;\Theta^{N-d} (\D (H_{d+1},\dots,H_{N})) \,\Theta^{d}(\D(H_{1},\dots,H_{d})) \\
& \le & \binom{N}{d} \Phi(K_m)^{-N} \int_{{\cal H}_{K_m}^{d}} [\Phi(K_m) - \Phi(K)- c_{24}\gamma  m ]^{N-d}\, \Theta^{d}(\D (H_{1},\dots,H_{d})) \\
& = & \binom{N}{d} \Phi(K_m)^{d-N} \left[\Phi(K_m) - \Phi(K)- c_{24}\gamma m\right]^{N-d}.
\end{eqnarray*}
With (\ref{1}) (for $\varepsilon=0$) and (\ref{2}) this gives
\begin{eqnarray*}
& & q_0(m)\\
& & \le \sum_{N=d+1}^{\infty} \frac{\Phi(K_m)^{N}}{N!} \exp\left[-\Phi(K_m)\right] \binom{N}{d} \Phi(K_m)^{d-N} \left[\Phi(K_m) - \Phi(K)-c_{24}\gamma m\right]^{N-d} \\
& & = \frac{1}{d!} \Phi(K_m)^{d} \exp[-\Phi(K_m)] \sum_{N=d+1}^{\infty} \frac{1}{(N-d)!} \left[\Phi(K_m)- \Phi(K)- c_{24}\gamma m \right]^{N-d} \\
& & \le \frac{1}{d!} \Phi(K_m)^{d} \exp\left[- \Phi(K)-c_{24}\gamma m\right] \\
& & = \frac{1}{d!} \left(2\gamma (m+1) \int_{\Sd} h(K,\ur)\,\varphi(\D\ur) \right)^d \exp\left[- \Phi(K)-c_{24}\gamma m\right]\\
& & \le c_{22}  \exp\left[- \Phi(K)-c_{23} \gamma m\right]
\end{eqnarray*}
with $c_{23}=c_{24}/2$, say. \qed

\vspace{3mm}

\noindent{\bf Lemma 5.} {\em Let $0<\varepsilon\le 1$. Then, for $m\in\mathbb{N}$,}
$$ q_\varepsilon(m) \le c_{25}(\gamma m)^d\exp\left[-\Phi(K)- 2\gamma \mu(K,\varphi,\varepsilon)\right].$$

\vspace{3mm}

\noindent{\em Proof.}
With $H_1,\dots,H_N \in {\cal H}_{K_m}$ and $P = P(H_{(N)})\in {\cal K}^d_0(m)$ as in the previous proof, the polytope $P$ has a vertex $\x\in K_m\setminus K(\varepsilon)$. This vertex is the intersection of $d$ facets of $P$. Hence, there exists an index set $J \subset \{1,\dots,N\}$ with $d$ elements such that
\[ \{\x\} = \bigcap_{j \in J} H_j. \]
There exists a point $\y\in {\rm bd}\,K(\varepsilon)$ such that
$$ \Phi({\rm conv}(K\cup\{\x\})) \ge \Phi({\rm conv}(K\cup\{\y\})) = \Phi(K^{\ys}) \ge \Phi(K) + 2\gamma \mu(K,\varphi,\varepsilon),$$
where the last inequality follows from (\ref{n2}) and (\ref{n4}), together with the monotonicity of $\Phi$. This gives
\begin{eqnarray*}
\int_{{\cal H}_{K_m}} {\bf 1} \{H\cap{\rm conv}(K\cup\{\x\})=\emptyset\}\,\Theta(\D H) &=& \Phi(K_m)- \Phi({\rm conv}(K\cup\{\x\}))\\
&\le& \Phi(K_m)- \Phi(K) - 2\gamma \mu(K,\varphi,\varepsilon).
\end{eqnarray*}
We write $\x=\x(H_1,\dots,H_d)$ for the intersection point of the hyperplanes $H_1,\dots,H_d$ (supposed in general position) and obtain
\begin{eqnarray*}
p(N,m,\varepsilon) &\le & \binom{N}{d} \Phi(K_m)^{-N} \int_{{\cal H}_{K_m}^d}  {\bf 1} \{ \x(H_1,\dots,H_d) \in K_m\setminus K(\varepsilon)\}\\
& & \int_{{\cal H}_{K_m}^{N-d}} {\bf 1}\{ H_i \cap {\rm conv}(K \cup \{\x(H_1,\dots,H_d)\}) = \emptyset \mbox{ for } i=d+1,\dots,N \}\\
& & \times \; \Theta^{N-d}(\D(H_{d+1},\dots,H_N))\, \Theta^d(\D(H_1\dots,H_d))\\
& \le & \binom{N}{d}\Phi(K_m)^{d-N}\left[\Phi(K_m) - \Phi(K)- 2\gamma \mu(K,\varphi,\varepsilon) \right]^{N-d}.
\end{eqnarray*}
Similarly as in the proof of Lemma 4, summation over $N$ gives
\begin{eqnarray*}
& & q_\varepsilon(m)\\
& & \le \sum_{N=d+1}^{\infty} \frac{\Phi(K_m)^{N}}{N!} \exp\left[-\Phi(K_m)\right] \binom{N}{d} \Phi(K_m)^{d-N} \left[\Phi(K_m) - \Phi(K)-2\gamma \mu(K,\varphi,\varepsilon)\right]^{N-d} \\
& & \le \frac{1}{d!}\Phi(K_m)^d \exp \left[-\Phi(K)- 2\gamma \mu(K,\varphi,\varepsilon)\right]\\
& & \le c_{25} (\gamma m)^d\exp \left[-\Phi(K)- 2\gamma \mu(K,\varphi,\varepsilon)\right].
\end{eqnarray*} \qed

\vspace{3mm}

\noindent{\em Proof of Theorem $1$.} We have
\begin{eqnarray*} 
\rP\left\{\delta(K,Z_K)>\varepsilon\right\} &=& \rP \left\{ Z_0 \not\subset K(\varepsilon) \mid K \subset Z_0\right\}\\
&=& \frac{\rP \left\{K \subset Z_0,\; Z_0 \not\subset K(\varepsilon) \right\}}{\rP\left\{K \subset Z_0\right\}}= \frac{\sum_{m=1}^\infty q_{\varepsilon}(m)}{\exp\left[-\Phi(K)\right]}.
\end{eqnarray*}
To estimate the last numerator, we choose $m_0$ according to Lemma 4 and use  Lemma 5 for $m\le m_0$ and Lemma 4 together with $q_{\varepsilon}(m) \le q_0(m)$ for $m>m_0$. By the assumptions of Theorem 1, relation $(\ref{n1})$ is satisfied. We obtain
$$ \rP \left\{ Z_0 \not\subset K(\varepsilon) \mid K \subset Z_0\right\} \le 
\sum_{m=1}^{m_0} c_{25}(\gamma m)^d \exp[-2\gamma \mu(K,\varphi,\varepsilon)] + \sum_{m>m_0} c_{22}\exp[-c_{23}\gamma m]. $$

The first sum can be estimated by
\begin{eqnarray} \label{4.1} 
& & \sum_{m=1}^{m_0} c_{25}(\gamma m)^d \exp[-2\gamma \mu(K,\varphi,\varepsilon)] \\
& & \le c_{25} m_0^{d+1} \gamma^d \exp\left[ -\gamma \mu(K,\varphi,\varepsilon)\right]\exp\left[ -\gamma 
\mu(K,\varphi,\varepsilon)\right]\nonumber\\
& & \le c_{26}(\varepsilon)\exp\left[-\gamma \mu(K,\varphi,\varepsilon)\right],\nonumber
\end{eqnarray}
since $\mu(K,\varphi,\varepsilon)>0$ by condition (\ref{n1}) and Lemma 1. 

The second sum can be estimated by
$$
\sum_{m>m_0} c_{22}\exp[-c_{23}\gamma m] \le c_{22}\exp[-c_{23}\gamma]\sum_{m>m_0} \exp[-c_{23}(m-1)]
\le c_{27}\exp\left[-c_{23}\gamma\right],
$$
where we have used that $\gamma\ge 1$ (by assumption) and that the last sum converges. Both estimates together yield (\ref{5.0}).
\qed

\section{Proofs of Theorems 2 and 3}

Under the assumptions (\ref{n7}) or (\ref{n3}), we can conclude from Lemma 1 that $\mu(K,\varphi,\varepsilon) \ge c_{28} \varepsilon^\alpha$ with suitable $\alpha\le d$. Therefore, in estimating (\ref{4.1})  we can use that
$$ \gamma^d \exp\left[ -\gamma \mu(K,\varphi,\varepsilon)\right] \le \gamma^d \exp\left( -\gamma c_{28}\varepsilon^\alpha\right)\le c_{29}\varepsilon^{-d\alpha}.$$
This gives
$$ \sum_{m=1}^{m_0} c_{25}(\gamma m)^d \exp[-2\gamma \mu(K,\varphi,\varepsilon)] \le  c_{30}\varepsilon^{-d\alpha}\exp\left(-c_{31}\gamma \varepsilon^\alpha\right).$$
The estimation of the second sum above remains unchanged. Hence, under the assumptions of Theorem 2 and with $\gamma=n$, we can conclude that
$$ \rP\left\{\delta(K,Z_K^{(n)})>\varepsilon\right\} \le c_{32} \varepsilon^{-d\alpha} \exp \left(-c_{33}n \varepsilon^\alpha \right).$$
We choose
$$ C>\frac{d+1}{c_{33}}$$
and put
$$ \varepsilon_n:= \left(\frac{C\log n}{n}\right)^{1/\alpha}.$$
Then
\begin{eqnarray}\label{n10}
\sum_{n=1}^\infty  \rP\left\{\delta(K,Z_K^{(n)})>\varepsilon_n\right\} &\le& \sum_{n=1}^\infty c_{32}\left(\frac{n}{C\log n} \right)^d \exp\left(-c_{33} C\log n\right) \nonumber\\
&=& c_{34}\sum_{n=1}^\infty (\log n)^{-d} n^{d-c_{33}C} <\infty.
\end{eqnarray}
The Borel--Cantelli lemma gives
$$ \rP \left\{\delta(K,Z_K^{(n)})>\varepsilon_n\mbox{ for infinitely many }n\right\}=0,$$
hence
$$ \rP \left\{\delta(K,Z_K^{(n)})\le\left(\frac{C\log n}{n}\right)^{1/\alpha} \mbox{ for sufficiently large }n\right\}=1.$$
This completes the proof of Theorem 2.\qed

\vspace{3mm}

\noindent{\em Proof of Theorem $3$.} 
Since $K$ slides freely in some ball, say of radius $R$, there is a convex body $L$ with $K+L=RB^d$ (\cite[Theorem 3.2.2]{Sch14}). From the polynomial expansion of $S_{d-1}(K+L,\cdot)$ (\cite[(5.18)]{Sch14}) it follows that $S_{d-1}(K,\cdot) \le S_{d-1}(RB^d,\cdot)=R^{d-1}\sigma$. Together with the assumption (\ref{n9}) this shows that (\ref{n3}) is satisfied. Therefore, Theorem 2 yields that 
\begin{equation}\label{5.0y}
\delta(K,Z_K^{(n)})={\rm O}\left(\left(\frac{\log n}{n}\right)^{\frac{2}{d+1}}\right)\quad\mbox{almost surely,}
\end{equation}
as $n\to\infty$. 

Let $0<\varepsilon < c_{20}$ (with $c_{20}$ as in Lemma 3). According to Lemma 3, we can choose
$$ m=m(\varepsilon) \ge c_{21}\varepsilon^{-1/2} $$
points $\x_1,\dots,\x_m\in{\rm bd}\,K(\varepsilon)$ such that the segment joining any two of them intersects the interior of $K$. Let $n\in{\mathbb N}$. Suppose that $\delta(K,Z_K^{(n)})<\varepsilon$. Then each point $\x_i$ is strictly separated from $K$ by some hyperplane from $X_n$. Let ${\mathcal A}_i\subset {\mathcal H}^d$ be the set of hyperplanes strictly separating $\x_i$ and $K$. By the choice of the points $\x_1,\dots,\x_m$, the sets ${\mathcal A}_1,\dots,{\mathcal A}_m$ are pairwise disjoint. Since $X_n$ is a Poisson process, the processes $X_n\fed {\mathcal A}_1,\dots,X_n\fed {\mathcal A}_m$ are stochastically independent (e.g., \cite[Theorem 3.2.2]{SW08}). It follows that
\begin{eqnarray*}
\rP \{ \delta(K,Z_K^{(n)}) < \varepsilon\} &\le& \rP \{ X_n({\mathcal A}_i)\ge 1\mbox{ for }i=1,\dots,m\}\\
&=& \prod_{i=1}^m \rP \{ X_n({\mathcal A}_i)\ge 1\}= \prod_{i=1}^m \left[ 1-\rP \{ X_n({\mathcal A}_i) =0\}\right]\\
&=& \prod_{i=1}^m \left( 1-\exp[-\Theta_n({\mathcal A}_i)] \right),
\end{eqnarray*}
where $\Theta_n$ is the intensity measure of $X_n$. Since the assumptions on $K$ in Lemma 2 are satisfied, we can conclude that
\begin{eqnarray*}
\Theta_n({\mathcal A}_i) &=& \Theta_n({\mathcal H}_{K^{\xs_i}}) - \Theta_n({\mathcal H}_K)\\
&=& 2n\int_{{\mathbb S}^{d-1}} [h(K^{\xs_i},\ur)-h(K,\ur)]\,\varphi(\D \ur)\\
&\le& 2nc_{17}\varepsilon^{(d+1)/2}.
\end{eqnarray*}
This gives
$$ \rP \left\{ \delta(K,Z_K^{(n)}) < \varepsilon\right \} \le \left[1-\exp\left(-2c_{17}n \varepsilon^{(d+1)/2}\right) \right]^{m(\varepsilon)}.$$
Now we choose
$$ \varepsilon_n^{(d+1)/2} = \frac{c\log n}{n}$$
with
$$ 0<c<\frac{1}{4c_{17}(d+1)}.$$ 
Then
$$ \rP \left\{ \delta(K,Z_K^{(n)}) < \varepsilon_n\right \} \le \left(1- n^{-2c_{17}c} \right)^{m(\varepsilon_n)}$$
with
$$ m(\varepsilon_n) \ge c_{21}\varepsilon_n^{-1/2} = c_{21}\left(\frac{n}{c\log n}\right)^{1/(d+1)}> c_{35} n^{1/(2d+2)}$$
for sufficiently large $n$. With $p:= 2c_{17}c$ and $q:= 1/(2d+2)$ we have $q>p$ and
$$ \left(1-n^{-2c_{17}c}\right)^{m(\varepsilon_n)} <\left(1-\frac{1}{n^p}\right)^{c_{35}n^q} =\left[\left(1-\frac{1}{n^p}\right)^{n^p\cdot n^{q-p}}\right]^{c_{35}} \le ({\rm e}^{-c_{35}})^{n^{q-p}}.$$
It follows that
$$ \sum_{n=1}^\infty \rP \left\{\delta(K,Z_K^{(n)}) <\left(\frac{c\log n}{n}\right)^{\frac{2}{d+1}}\right\} <\infty.$$
From the Borel--Cantelli lemma we conclude that
$$ \rP \left\{\delta(K,Z_K^{(n)}) <\left(\frac{c\log n}{n}\right)^{\frac{2}{d+1}}\mbox{ for infinitely many }n\right\} =0$$
and hence 
$$ \rP \left\{\delta(K,Z_K^{(n)})  \ge \left( \frac{c\log n}{n}\right)^{\frac{2}{d+1}}
\mbox{ for almost all }n \right\} =1.$$
Together with (\ref{5.0y}), this completes the proof of Theorem 3.\qed

Finally, we construct a directional distribution exhibiting the property (\ref{C1}) for a given convex body $K$. We do that at this stage, since arguments appearing in the previous proofs are employed.

As explained before (\ref{C1}), we assume that a decreasing sequence $(\varepsilon_n)_{n\in{\mathbb N}}$ with $\lim_{n\to\infty}\varepsilon_n=0$ is given. For $n\in{\mathbb N}$, let $X_n$ be a stationary Poisson hyperplane process with intensity $n$ and directional distribution $\varphi$, to be constructed. 

The $d$-dimensional convex body $K$ contains some ball touching the boundary, hence there are a number $r>0$ and a point $\x\in{\rm bd}\,K$ such that $\x$ is contained in a ball of radius $r$ that is contained in $K$. Let $N(\x)$ be the unique outer unit normal vector of $K$ at $\x$ and let $\y=\x+\varepsilon N(\x)$. Let $n\in{\mathbb N}$, and suppose that $\delta(K,Z_K^{(n)})<\varepsilon_n$. Then the point $\y$ is strictly separated from $K$ by some hyperplane of $X_n$. Similarly as in the proof of Theorem 3, this yields
$$ \rP\left\{\delta(K,Z_K^{(n)})<\varepsilon_n\right\}\le 1-\exp[-2n\varepsilon_n\varphi(S_n)]$$
with
$$ S_n:=\left\{\ur\in\Sd:\langle \ur,N(\x) \rangle\ge \frac{r}{r+\varepsilon_n}\right\}.$$

It is easy to construct an even positive measurable function $g$ on $\Sd$ such that the measure $\varphi$ defined by $\D\varphi=g \,\D\sigma$ is a probability measure and that 
$$
2n\varepsilon_n\varphi(S_n) < |\log(1-n^{-2})|
$$
for all $n\in{\mathbb N}$ (for example, $g$ can be a suitable constant on $S_n\setminus S_{n+1}$). The directional distribution $\varphi$ then satisfies
$$ 1-\exp[-2n\varepsilon_n\varphi(S_n)]<\frac{1}{n^2}$$ and hence
$$ \sum_{n=1}^\infty \rP\left\{\delta(K,Z_K^{(n)})<\varepsilon_n\right\}<\infty.$$
As in the proof of Theorem 3, this yields (\ref{C1}).

\noindent Authors' addresses:\\[2mm]
Daniel Hug\\
Karlsruhe Institute of Technology, Department of Mathematics\\
D-76128 Karlsruhe, Germany\\
E-mail: daniel.hug@kit.edu\\[3mm]
Rolf Schneider\\
Mathematisches Institut, Albert-Ludwigs-Universit{\"a}t Freiburg\\
D-79104 Freiburg i. Br., Germany\\
E-mail: rolf.schneider@math.uni-freiburg.de

\end{document}